\documentclass[12pt]{amsart}
\usepackage[T1]{fontenc}
\usepackage[utf8]{inputenc}
\usepackage{lmodern}
\usepackage{amsmath,amssymb,mathtools}
\usepackage{microtype}
\usepackage{enumitem}
\usepackage{xcolor}

\usepackage[hidelinks]{hyperref}
\numberwithin{equation}{section}
\newtheorem{theorem}{Theorem}[section]
\newtheorem{proposition}[theorem]{Proposition}
\newtheorem{lemma}[theorem]{Lemma}
\newtheorem{corollary}[theorem]{Corollary}
\theoremstyle{remark}

\newcommand{\HH}{\mathcal H}
\newcommand{\BB}{\mathcal B}
\newcommand{\CC}{\mathbb C}
\newcommand{\MM}{\mathcal M}
\newcommand{\SSS}{\mathbb S^2}
\newcommand{\norm}[1]{\left\|#1\right\|}
\DeclareMathOperator{\diag}{diag}
\DeclareMathOperator{\rank}{rank}
\DeclareMathOperator{\ran}{ran}
\allowdisplaybreaks[1]
\title[Similarity cancellation and primitive roots]{Similarity cancellation and primitive operator roots}
\author{Mohamed Amine Aouichaoui}
\address{Department of Mathematics, Faculty of Sciences of Monastir, University of Monastir, Tunisia}
\email{aouichaoui.mohamedamine@fsm.u-monastir.tn}
\email{amine.aouichaoui@fsm.rnu.tn}
\subjclass[2020]{Primary 47A45, 47A65; Secondary 46L10, 47B10}
\keywords{Similarity, finite amplification, primitive operator roots, compact operators, von Neumann algebras}
\begin{document}
\begin{abstract}
Kaplansky's second test problem for similarity asks whether $A\oplus A\sim B\oplus B$ implies $A\sim B$ for bounded operators on a separable infinite-dimensional Hilbert space. We answer this question in the negative, even for compact operators whose singular values satisfy any prescribed positive upper bounds tending to zero. We also show that cancellation can fail under the stronger hypothesis $J_2(T)\sim J_2(S)$ considered by Marcoux, Radjavi, Troscheit, and Zhang \cite{MRTZ}. For these operators, finite amplifications and primitive roots are similar exactly at even multiplicities and even orders, respectively. In contrast, cancellation in fixed matrix size admits a uniform bound on the condition number. This extends the finite type~$I$ results of Jiang, Ma, Shi, and Zhang \cite{JMSZ}. The bound also holds for simultaneous similarity of arbitrary families of matrices and can be chosen to tend to one as the amplified bound tends to one. We conclude that, in a finite type~$I$ von Neumann algebra, cancellation of doubles and cancellation of primitive square roots are each equivalent to boundedness of the orders of its nonzero homogeneous summands.

\end{abstract}
\maketitle

\section{Introduction}\label{sec:introduction}

Let $H$ be a complex separable infinite-dimensional Hilbert space. We denote by $\BB(H)$ the algebra of bounded operators on $H$. Two operators $A$ and $B$ are similar, written $A\sim B$, if $YA=BY$ for some bounded invertible operator $Y$. We use the same notation for operators on different Hilbert spaces, with $Y$ acting between those spaces. All Hilbert space direct sums are orthogonal. For a positive integer $q$, write
\[
A^{(q)}=\underbrace{A\oplus\cdots\oplus A}_{q\text{ copies}}.
\]

In his study of abelian groups, Kaplansky asked whether $G\oplus G\cong F\oplus F$ forces $G\cong F$ \cite{Kaplansky}. This is the second of his three test problems for direct sums. Kadison and Singer \cite{KS} brought these questions to Hilbert space operators, and Azoff \cite{Azoff} studied them in several categories of operator algebras. {Related versions of Kaplansky's test problems for quasisimilarity were studied by Bercovici \cite{Bercovici} and, for bitriangular operators, by Davidson and Herrero \cite{DH}.} The answer depends on the equivalence relation being considered. For unitary equivalence, Kadison and Singer proved that
\[
A\oplus A\simeq_u B\oplus B\quad\Longrightarrow\quad A\simeq_u B,
\]
where $\simeq_u$ denotes unitary equivalence. There is also a multiplicity theory for approximate unitary equivalence. Hadwin's rank criterion for unital representations on a separable Hilbert space \cite[Theorem~2.5]{Hadwin} gives a formulation of Voiculescu's theorem suited to this question; see also \cite[\S II.5]{Davidson}. In particular, finite amplifications cancel for approximate unitary equivalence in $\BB(H)$, as shown in \cite[Proposition~4.4]{MRTZ}. Related questions concerning distances between unitary orbits of matrix amplifications were studied by Marcoux, Sarkowicz, and Zhang \cite{MSZ}. Their results concern a different quantity from the condition number considered here. For similarity, the question is
\begin{equation}\label{eq:kaplansky}
A\oplus A\sim B\oplus B\quad\Longrightarrow\quad A\sim B\;?
\end{equation}
{
Jiang \cite[Theorem~4.6]{Jiang} proved cancellation of finite amplifications for strongly irreducible Cowen--Douglas operators. Jiang, Guo, and Ji \cite{JGJ} developed the similarity classification of Cowen--Douglas operators using the ordered $K$-group of the commutant algebra. For weighted Hardy spaces of polynomial growth, Hou and Jiang \cite[Theorem~B]{HJ} proved cancellation of doubles for the Hermitian holomorphic vector bundles induced by functions analytic in a neighbourhood of the closed unit disc.
\par}

Marcoux, Radjavi, Troscheit, and Zhang record this question in \cite[\S1.3]{MRTZ}, including the case in which $A$ and $B$ are compact.

The same authors consider a stronger hypothesis. Following their notation, for $n\ge2$ let
\begin{equation}\label{eq:rootdef}
J_n(A)(h_1,\ldots,h_n)=(h_2,\ldots,h_n,Ah_1).
\end{equation}
They call $J_n(A)$ the primitive $n$th root of $A^{(n)}$. In particular,
\[
J_2(A)=\begin{pmatrix}0&I\\A&0\end{pmatrix},
\qquad J_n(A)^n=A^{(n)}.
\]
Thus, similarity of $J_2(A)$ and $J_2(B)$ implies similarity of the doubles. The question in \cite[Remark~5.3]{MRTZ} is whether
\begin{equation}\label{eq:rootquestion}
J_2(A)\sim J_2(B)\quad\Longrightarrow\quad A\sim B.
\end{equation}
Their Theorem~4.10 proves this implication when one of the operators is compact. A counterexample to cancellation in a unital ring, due to J.~Bell, is recorded in \cite[Example~5.4]{MRTZ}. That algebraic example does not settle either question in $\BB(H)$.

We answer both questions in the negative. Our counterexample to \eqref{eq:rootquestion} is noncompact, as required by the theorem just cited. The same pair distinguishes even from odd multiplicities, both for direct sums and for primitive roots.

\begin{theorem}\label{thm:bounded}
For every $\varepsilon>0$, there are noncompact operators $T,S\in\BB(H)$ with $\norm{T},\norm{S}<\varepsilon$ such that
\begin{align}
T^{(q)}\sim S^{(q)}&\quad\Longleftrightarrow\quad q\text{ is even}
&& (q\ge1),\label{eq:parity}\\
J_n(T)\sim J_n(S)&\quad\Longleftrightarrow\quad n\text{ is even}
&& (n\ge2).\label{eq:rootparity}
\end{align}
The operators have the same point spectrum. Every eigenspace is one-dimensional, and the eigenvectors of each operator have dense linear span.
\end{theorem}

The answer to \eqref{eq:kaplansky} is negative even when both operators are compact. In fact, their singular values can lie below any prescribed positive sequence tending to zero. For a compact operator $R$, let $s_j(R)$ denote its singular values in nonincreasing order, counted with multiplicity.

\begin{theorem}\label{thm:compact}
Let $(a_j)$ be a nonincreasing sequence of positive numbers tending to zero. There are compact operators $K,L\in\BB(H)$ such that
\begin{equation}\label{eq:singularbounds}
s_j(K),s_j(L)\le a_j\qquad(j\ge1)
\end{equation}
and
\begin{equation}\label{eq:compactparity}
K^{(q)}\sim L^{(q)}\quad\Longleftrightarrow\quad q\text{ is even}
\qquad(q\ge1).
\end{equation}
Both operators are injective and have dense range. Their common nonzero spectrum consists of positive algebraically simple eigenvalues, and their eigenvectors have dense linear span.
\end{theorem}

For an arbitrary positive null sequence $(a_j)$, apply Theorem~\ref{thm:compact} to the nonincreasing sequence $b_j=\min_{1\le i\le j}a_i$. Since $0<b_j\le a_j$ and $b_j\to0$, the same upper bounds follow without a monotonicity assumption. In particular, $K$ and $L$ may belong to every Schatten class. Their positive eigenvalues should not be confused with positivity of the operators.

The counterexamples are orthogonal sums of diagonalizable matrices. In finite dimension, the Jordan form gives cancellation of finite amplifications; see \cite[\S1.3]{MRTZ}. For an infinite orthogonal sum, however, the similarities chosen on the summands must have uniformly bounded norms and inverse norms. Pointwise similarity without such bounds is not enough, as the examples of Deckard and Pearcy already show \cite[\S2]{DP}. We therefore also consider the condition number of a similarity. For matrices $A,B\in M_n(\CC)$, set
\begin{equation}\label{eq:similarity-number}
s(A,B)=\inf\{\norm{Y}\norm{Y^{-1}}:Y\in GL_n(\CC),\ YA=BY\},
\end{equation}
with $s(A,B)=\infty$ when $A$ and $B$ are not similar. For $n,q\ge1$ and $t\ge1$, put
\begin{equation}\label{eq:theta}
\theta_{n,q}(t)=\sup\{s(A,B):A,B\in M_n(\CC),\
 s(A^{(q)},B^{(q)})\le t\}.
\end{equation}
For $q=2$, this is the function introduced by Jiang, Ma, Shi, and Zhang \cite[Definition~4.1]{JMSZ}.

\begin{theorem}\label{thm:uniform}
For every $n,q\ge1$ and $t\ge1$,
\[
\theta_{n,q}(t)<\infty.
\]
\end{theorem}

Jiang, Ma, Shi, and Zhang relate this bound to similarity cancellation in homogeneous type~$\mathrm I_n$ von Neumann algebras \cite[Theorem~4.3]{JMSZ}. They prove cancellation under property~$(J)$ \cite[Theorem~3.27]{JMSZ} and, without this assumption, for $n\le3$ \cite[Proposition~4.7]{JMSZ}. Theorem~\ref{thm:uniform}, together with the continuous selection theorem of Deckard and Pearcy \cite{DP}, removes the extra assumption for every finite $n$ and every multiplicity $q$.

Two further conclusions follow from the proof. The same uniform bound, with a possibly different constant, holds when one operator is required to intertwine an arbitrary family of matrices simultaneously. The constant does not depend on the number of matrices in the family. For a single pair, we also prove
\begin{equation}\label{eq:near-one-intro}
\lim_{t\downarrow1}\theta_{n,q}(t)=1
\qquad(n,q\text{ fixed}).
\end{equation}
Thus, a similarity with condition number close to one remains so after cancellation, provided the matrix size is fixed. Proposition~\ref{prop:nearone} proves the same assertion for a common intertwiner of arbitrary families, uniformly in the index set.

The dependence on matrix size cannot be omitted. The finite-dimensional operators used in Theorems~\ref{thm:bounded} and~\ref{thm:compact} have uniformly controlled similarities after doubling, while every similarity at an odd multiplicity has arbitrarily large condition number. Consequently, there is a fixed $t_0$ such that
\begin{equation}\label{eq:size-dependence}
\theta_{n,2}(t_0)<\infty\quad\text{for every }n,
\qquad
\sup_{n\ge1}\theta_{n,2}(t_0)=\infty.
\end{equation}
These pairs also show that $s(A^{(2)},B^{(2)})<s(A,B)$ can occur, contrary to the equality suggested in \cite[Remark~4.2]{JMSZ}; see Corollary~\ref{prop:finite-bad}.
The relation between a dimension-independent bound and cancellation in finite type~I algebras was established in \cite[Corollary~4.6]{JMSZ}. We decide cancellation for each bounded product $\prod_jM_{d_j}(\CC)$: both the double and primitive-square-root questions have affirmative answers for all pairs in the product if and only if $\sup_jd_j<\infty$. More generally, for a finite type~I von Neumann algebra, the same two properties are equivalent to a common finite bound on the orders of its nonzero homogeneous summands. The centres of those summands may be diffuse.

The proofs of the positive and negative results both concern uniform bounds for finite-dimensional intertwiners. To prove Theorem~\ref{thm:uniform}, we take a limit of the commutants of $A_k\oplus B_k$ inside a fixed matrix space. The limit is a finite-dimensional algebra, so Krull--Schmidt cancellation applies. The resulting intertwiner can be approximated by exact intertwiners for the original pairs.

For the counterexamples, the starting point is the even and odd decomposition of a Sobolev space on the sphere. Restriction to a finite set gives two Hilbert norms on the same vector space of functions. An intertwiner between multiplication operators with distinct eigenvalues must act pointwise. Applying it to constant functions and lifting to the odd Sobolev space produces a continuous odd matrix-valued function. At an odd multiplicity its determinant must vanish, by Borsuk--Ulam. Compactness of the inclusion into $C(\SSS)$ then forces the condition number to be large. At multiplicity two, a fixed invertible block operator induces similarities on all the quotients. Its entries are chosen to satisfy the additional identity for primitive square roots. Taking orthogonal sums with disjoint finite spectra gives the two pairs of operators.

The paper is organised as follows. Section~\ref{sec:fixed} proves the fixed-size bound, its simultaneous version, and \eqref{eq:near-one-intro}, and then treats homogeneous type~$\mathrm I_n$ algebras. Section~\ref{sec:quotients} constructs the Hilbert space quotients and their intertwiners. Section~\ref{sec:examples} proves Theorems~\ref{thm:bounded} and~\ref{thm:compact} by taking orthogonal sums. Section~\ref{sec:products} compares bounded and unbounded matrix sizes and characterises cancellation in finite type~I von Neumann algebras.

\section{Cancellation in fixed matrix size}\label{sec:fixed}

We first prove the uniform estimate for matrices. The argument uses cancellation inside a finite-dimensional algebra and compactness of the set of orthogonal projections on a fixed matrix space.

By Krull--Schmidt cancellation, $P^{\oplus q}\cong Q^{\oplus q}$ implies $P\cong Q$ for finite-dimensional modules over a finite-dimensional algebra; see \cite[Theorems~4.2 and~5.5]{Krause}.

We apply this result to right ideals in finite-dimensional operator algebras, without assuming self-adjointness. Let $\mathcal C$ be a finite-dimensional unital algebra and let $e,f\in\mathcal C$ be idempotents. Suppose
\begin{equation}\label{eq:stable-corners}
U\in M_q(f\mathcal C e),\qquad V\in M_q(e\mathcal C f),
\qquad VU=e^{(q)},\quad UV=f^{(q)}.
\end{equation}
Left multiplication by $U$ identifies $(e\mathcal C)^{\oplus q}$ with $(f\mathcal C)^{\oplus q}$. Cancellation gives $e\mathcal C\cong f\mathcal C$. An isomorphism between these right modules is left multiplication by an element $u\in f\mathcal C e$, and its inverse is left multiplication by some $v\in e\mathcal C f$. Thus,
\begin{equation}\label{eq:single-corners}
vu=e,\qquad uv=f.
\end{equation}
We are now ready to prove Theorem~\ref{thm:uniform}.

\begin{proof}[Proof of Theorem~\ref{thm:uniform}]
Fix $n,q$ and $t$, and suppose that $\theta_{n,q}(t)=\infty$. For each $k$, choose $A_k,B_k\in M_n(\CC)$ with
\begin{equation}\label{eq:bad-pairs}
s(A_k^{(q)},B_k^{(q)})\le t,\qquad s(A_k,B_k)>k.
\end{equation}
By the definition of the infimum, there is an invertible $X_k\in M_{nq}(\CC)$ such that
\begin{equation}\label{eq:amplified-equation}
X_k A_k^{(q)}=B_k^{(q)}X_k,
\qquad \norm{X_k}\norm{X_k^{-1}}<t+1.
\end{equation}
Multiplying $X_k$ by a positive scalar, we may assume
\begin{equation}\label{eq:uniform-X}
\norm{X_k}=\norm{X_k^{-1}}<\sqrt{t+1}.
\end{equation}

Consider the commutant
\[
\mathcal C_k=\{A_k\oplus B_k\}'\subset M_{2n}(\CC)
\]
and the two fixed projections
\[
e=\begin{pmatrix}I_n&0\\0&0\end{pmatrix},
\qquad
f=\begin{pmatrix}0&0\\0&I_n\end{pmatrix}.
\]
Both projections belong to every $\mathcal C_k$. Equip $M_{2n}(\CC)$ with the Hilbert--Schmidt inner product, and let $P_k$ be the orthogonal projection onto the linear subspace $\mathcal C_k$. The maps $P_k$ act on a fixed finite-dimensional Hilbert space and have norm one. After passing to a subsequence, we may suppose that
\[
P_k\longrightarrow P,\qquad X_k\longrightarrow X
\]
in norm. The limit $P$ is an orthogonal projection. For every $h\in\CC^{nq}$, \eqref{eq:uniform-X} gives $\norm{X_kh}\ge(t+1)^{-1/2}\norm{h}$. The same lower bound holds for $X$. Since $X$ is square and finite-dimensional, it is invertible, and $X_k^{-1}\to X^{-1}$.

Put $\mathcal C=\ran P$. We first check that $\mathcal C$ is an algebra. If $R,S\in\mathcal C$, then $R_k=P_kR\to R$ and $S_k=P_kS\to S$. Since $R_kS_k\in\mathcal C_k$,
\[
P_k(R_kS_k)=R_kS_k.
\]
Taking limits gives $P(RS)=RS$. Hence, $RS\in\mathcal C$. Moreover, $I,e,f\in\mathcal C$, since each of these matrices is fixed by every $P_k$. Thus, $\mathcal C$ is a finite-dimensional unital algebra. 

Write $X_k=[X_{ij}^{(k)}]$ and $X_k^{-1}=[Z_{ij}^{(k)}]$ in $q\times q$ blocks of size $n$. Equation~\eqref{eq:amplified-equation} and its inverse give
\[
\begin{pmatrix}0&0\\X_{ij}^{(k)}&0\end{pmatrix}\in f\mathcal C_k e,
\qquad
\begin{pmatrix}0&Z_{ij}^{(k)}\\0&0\end{pmatrix}\in e\mathcal C_k f.
\]
If $R_k\in\mathcal C_k$ and $R_k\to R$, then $P_kR_k=R_k$ implies $PR=R$. Therefore the corresponding corner matrices formed from the entries of $X$ and $X^{-1}$ belong to $f\mathcal C e$ and $e\mathcal C f$. Denote the resulting $q\times q$ matrices over $\mathcal C$ by $U$ and $V$. The equations $X^{-1}X=XX^{-1}=I_{nq}$ give \eqref{eq:stable-corners}.

By \eqref{eq:single-corners}, there are $u\in f\mathcal C e$ and $v\in e\mathcal C f$ with $vu=e$ and $uv=f$. Write
\[
u=\begin{pmatrix}0&0\\Y&0\end{pmatrix},
\qquad
v=\begin{pmatrix}0&Z\\0&0\end{pmatrix}.
\]
Then, $ZY=YZ=I_n$, so $Y$ is invertible.

For each $k$, the matrix $u_k=fP_k(u)e$ belongs to $f\mathcal C_k e$ and tends to $u$. Write its lower-left block as $Y_k$. We have
\[
Y_k A_k=B_kY_k,\qquad Y_k\longrightarrow Y.
\]
For all sufficiently large $k$, $Y_k$ is invertible and $Y_k^{-1}\to Y^{-1}$. Thus, $s(A_k,B_k)\le\norm{Y_k}\norm{Y_k^{-1}}$ stays bounded, contradicting \eqref{eq:bad-pairs}.
\end{proof}

The proof does not require the sequences $(A_k)$ and $(B_k)$ to converge, or even to be bounded. It is the dimension of the ambient commutant space that is fixed. Also, the limiting algebra $\mathcal C$ need not itself be a commutant.

\subsection{Simultaneous similarity}
The proof just given also applies to a joint commutant. This allows the same operator to intertwine any number of matrices, without a commutativity assumption on the family.

\begin{theorem}\label{thm:simultaneous}
Fix $n,q\ge1$ and $t\ge1$. There is a constant $C=C(n,q,t)$ with the following property. Let $(A_\alpha)_{\alpha\in I}$ and $(B_\alpha)_{\alpha\in I}$ be families in $M_n(\CC)$, indexed by the same nonempty set. Suppose that one invertible matrix $X$ satisfies
\[
XA_\alpha^{(q)}=B_\alpha^{(q)}X\quad(\alpha\in I),
\qquad \norm{X}\norm{X^{-1}}\le t.
\]
Then, there is an invertible $Y\in M_n(\CC)$ such that
\[
YA_\alpha=B_\alpha Y\quad(\alpha\in I),
\qquad \norm{Y}\norm{Y^{-1}}\le C.
\]
The constant is independent of $I$ and of the two families.
\end{theorem}

We apply the preceding commutant argument to all the equations at once.

\begin{proof}
Suppose that no such constant exists. For each $k$, choose two families $(A_{k,\alpha})_{\alpha\in I_k}$ and $(B_{k,\alpha})_{\alpha\in I_k}$ with a common amplified intertwiner $X_k$ of condition number at most $t$, but no common intertwiner of condition number at most $k$. After balancing $X_k$, we have $\norm{X_k}=\norm{X_k^{-1}}\le\sqrt t$.

Apply the proof of Theorem~\ref{thm:uniform} with the joint commutants
\[
\mathcal C_k=
\{A_{k,\alpha}\oplus B_{k,\alpha}:\alpha\in I_k\}'
\subset M_{2n}(\CC).
\]
Each block of $X_k$ satisfies every intertwining equation at the $k$th stage, and the blocks of $X_k^{-1}$ satisfy the reverse equations. Thus the same limiting algebra and corner argument give an element $u\in f\mathcal C e$ with invertible lower-left block. The lower-left blocks $Y_k$ of $fP_k(u)e$ converge to that block and satisfy
\[
Y_kA_{k,\alpha}=B_{k,\alpha}Y_k
\qquad(\alpha\in I_k).
\]
For all sufficiently large $k$, the matrices $Y_k$ are invertible and their condition numbers are bounded, contrary to the choice of the families.
\end{proof}

The index set may be infinite. The proof takes place in $M_{2n}(\CC)$ and uses only the common intertwining equations, not a bound on the number of those equations.

\subsection{Condition numbers close to one}
For an invertible matrix, condition number one means that a scalar multiple is unitary. The next proposition shows that this remains a uniform limiting statement after cancellation in a fixed size.

\begin{proposition}\label{prop:nearone}
For every fixed $n,q\ge1$,
\begin{equation}\label{eq:nearone}
\theta_{n,q}(1)=1,
\qquad
\lim_{t\downarrow1}\theta_{n,q}(t)=1.
\end{equation}
More generally, for every $\varepsilon>0$ there is $\delta>0$, depending only on $n,q,\varepsilon$, with the following property. Let $(A_\alpha)_{\alpha\in I}$ and $(B_\alpha)_{\alpha\in I}$ be families in $M_n(\CC)$ indexed by the same nonempty set. If an invertible $X\in M_{nq}(\CC)$ satisfies
\[
XA_\alpha^{(q)}=B_\alpha^{(q)}X\quad(\alpha\in I),
\qquad \norm{X}\norm{X^{-1}}\le1+\delta,
\]
then there is an invertible $Y\in M_n(\CC)$ such that
\[
YA_\alpha=B_\alpha Y\quad(\alpha\in I),
\qquad \norm{Y}\norm{Y^{-1}}\le1+\varepsilon.
\]
\end{proposition}

In this case the limiting algebra contains the adjoints of the corner entries. This allows us to choose a partial isometry before returning to the original commutants.

\begin{proof}
We prove the assertion for families first. If it fails, there are $\varepsilon>0$ and, for each $k$, families $(A_{k,\alpha})_{\alpha\in I_k}$ and $(B_{k,\alpha})_{\alpha\in I_k}$ with a common amplified intertwiner $X_k$ of condition number at most $1+k^{-1}$, but no common intertwiner of condition number at most $1+\varepsilon$. After balancing,
\[
\norm{X_k}=\norm{X_k^{-1}}\le\sqrt{1+k^{-1}}.
\]
Use the joint commutants from the proof of Theorem~\ref{thm:simultaneous} and the projections and subsequences from the proof of Theorem~\ref{thm:uniform}. We obtain $P_k\to P$, $X_k\to X$, and the unital algebra $\mathcal C=\ran P$. Here $X$ and $X^{-1}$ have norm at most one, so $X$ is unitary.

Let $\mathcal A=\mathcal C\cap\mathcal C^*$. This is a finite-dimensional unital $C^*$-algebra containing $e$ and $f$. The lower-left corner matrices formed from the blocks of $X$ belong to $\mathcal C$. Their adjoints also belong to $\mathcal C$, since $X^{-1}=X^*$ and the upper-right corner matrices formed from the blocks of $X_k^{-1}$ have limits in $\mathcal C$. Hence, the matrices $U,V$ of \eqref{eq:stable-corners} belong to the corresponding corners over $\mathcal A$, and $V=U^*$.

Apply \eqref{eq:single-corners} in $\mathcal A$. We obtain $u\in f\mathcal A e$ and $v\in e\mathcal A f$ with $vu=e$ and $uv=f$. Thus, $u$ is an invertible map from $e\CC^{2n}$ onto $f\CC^{2n}$. By continuous functional calculus in the corner $e\mathcal A e$ \cite[\S I.4]{Davidson}, its polar part
\[
w=u(u^*u)^{-1/2}
\]
belongs to $f\mathcal A e$ and satisfies $w^*w=e$ and $ww^*=f$. The inverse square root here is taken in $e\mathcal A e$, whose identity is $e$. Write
\[
w=\begin{pmatrix}0&0\\W&0\end{pmatrix};
\]
then $W$ is unitary.

The lower-left blocks $W_k$ of $fP_k(w)e$ satisfy
\[
W_kA_{k,\alpha}=B_{k,\alpha}W_k\qquad(\alpha\in I_k)
\]
and converge to $W$. For all sufficiently large $k$, they are invertible and
\[
\norm{W_k}\norm{W_k^{-1}}\longrightarrow1,
\]
contrary to the choice of the families.

Now fix $\varepsilon>0$ and choose $\delta>0$ as in the assertion for families. For a single pair with $s(A^{(q)},B^{(q)})\le1+\delta/2$, the definition of the infimum gives an amplified intertwiner with condition number less than $1+\delta$. The assertion just proved yields $s(A,B)\le1+\varepsilon$. Hence $\theta_{n,q}(1+\delta/2)\le1+\varepsilon$, which proves the limit in \eqref{eq:nearone}. Finally, $1\le\theta_{n,q}(1)\le\theta_{n,q}(t)$ for $t>1$, so $\theta_{n,q}(1)=1$.
\end{proof}

The matrix size is fixed throughout \eqref{eq:nearone}. Section~\ref{sec:products} shows that a bound valid in all sizes is a different question.

\subsection{Homogeneous type \texorpdfstring{$\mathrm I_n$}{I-n} algebras}

We next apply the fixed-size estimate to continuous matrix-valued functions. Deckard and Pearcy studied these algebras in \cite{DP63,DP}; the uniform norm requirement in their similarity theorem is the point needed here.

A homogeneous type~$\mathrm I_n$ von Neumann algebra is isomorphic to $M_n(\mathcal A)$ for an abelian von Neumann algebra $\mathcal A$. We may write $\mathcal A=C(\Omega)$, where $\Omega$ is compact and hyperstonean, hence Stonian; see \cite[\S6.6]{KR}. We use the following continuous selection result of Deckard and Pearcy \cite[Theorem~1, p.~858]{DP}: if $A,B\in M_n(C(\Omega))$ are pointwise similar on a dense subset of $\Omega$, with intertwiners whose norms and inverse norms are all less than one fixed number $M$, then there is an invertible $Y\in M_n(C(\Omega))$ with $YA=BY$ and
\[
\norm{Y},\norm{Y^{-1}}\le M.
\]

\begin{theorem}\label{thm:typeI}
Let $\MM$ be a homogeneous type~$\mathrm I_n$ von Neumann algebra, with $n<\infty$, and let $A,B\in\MM$. If $A^{(q)}\sim B^{(q)}$ in $M_q(\MM)$ for some positive integer $q$, then $A\sim B$ in $\MM$. In particular, $J_q(A)\sim J_q(B)$ in $M_q(\MM)$ implies $A\sim B$ for every $q\ge2$.
\end{theorem}

We can now apply the selection theorem to the pointwise similarities.

\begin{proof}
Write $\MM=M_n(C(\Omega))$ as above. Suppose that $X A^{(q)}=B^{(q)}X$ for an invertible $X\in M_q(\MM)$, and set $t=\norm{X}\norm{X^{-1}}$. Evaluation at $\omega\in\Omega$ gives
\[
s\bigl(A(\omega)^{(q)},B(\omega)^{(q)}\bigr)\le t.
\]
By Theorem~\ref{thm:uniform},
\[
s\bigl(A(\omega),B(\omega)\bigr)\le\theta_{n,q}(t)<\infty
\qquad(\omega\in\Omega).
\]
Choose $M>\sqrt{\theta_{n,q}(t)}$. For each $\omega$, the definition of the infimum gives a pointwise intertwiner with condition number less than $M^2$. A scalar multiple makes both its norm and its inverse norm less than $M$. The theorem of Deckard and Pearcy applies on all of $\Omega$ and gives a similarity in $\MM$. For primitive roots, take $q$th powers and apply the assertion just proved.
\end{proof}

\section{Hilbert quotients and finite-dimensional intertwiners}\label{sec:quotients}

We now construct finite-dimensional pairs with uniformly bounded similarities of their doubles, but without such a bound before amplification. We work with two fixed Hilbert spaces and take their quotients by the functions that vanish on a finite set.

\subsection{Even and odd quotient spaces}

We choose the Hilbert spaces so that evaluation is continuous and smooth multiplication is bounded.

Let $\HH=H^2(\SSS)$ be the complex Sobolev space on the unit sphere in $\mathbb R^3$, equipped with an orthogonally invariant Sobolev norm, normalised so that $\norm{1}_{\HH}=1$. Multiplication by every smooth function is bounded on $\HH$, and the natural inclusion
\begin{equation}\label{eq:compactinclusion}
\HH\longrightarrow C(\SSS)
\end{equation}
is compact. These are consequences of the smooth-multiplier, Sobolev embedding, and Rellich theorems; see \cite[Chapter~4, \S\S2--3]{Taylor}. In particular, compactness follows by factoring through $H^t(\SSS)$ for $1<t<2$, using Propositions~3.3--3.4 of that chapter. We use the continuous representative of each element of $\HH$.

The operator
\[
(\Gamma f)(x)=f(-x)
\]
is a unitary involution. Its orthogonal eigenspaces
\[
\HH_+=\ker(\Gamma-I),\qquad \HH_-=\ker(\Gamma+I)
\]
consist of the even and odd functions. A smooth even multiplier preserves both spaces, while a smooth odd multiplier maps each into the other.

Let $e_x^+$ and $e_x^-$ denote evaluation at $x$ on $\HH_+$ and $\HH_-$. Orthogonal invariance gives constants $a,b>0$ such that
\begin{equation}\label{eq:evalnorms}
\norm{e_x^+}=a,\qquad\norm{e_x^-}=b
\qquad(x\in\SSS).
\end{equation}
The first functional is nonzero because it evaluates the constant function. The second is nonzero because a smooth odd function can be chosen not to vanish at any prescribed point. By \eqref{eq:compactinclusion} and the Arzel\`a--Ascoli theorem \cite[Theorem~VI.3.8]{Conway}, the unit ball of $\HH_-$ is uniformly equicontinuous. Hence,
\begin{equation}\label{eq:evalcontinuous}
\sup_{|x-y|<\eta}\norm{e_x^- -e_y^-}\longrightarrow0
\qquad(\eta\downarrow0).
\end{equation}
The norm in this formula is the dual norm on $\HH_-$.

Let $E$ be a finite nonempty subset of $\SSS$ containing no antipodal pair. Let $H_E$ and $K_E$ be the restrictions to $E$ of $\HH_+$ and $\HH_-$, with the quotient norms
\begin{equation}\label{eq:quotientnorms}
\begin{split}
\norm{u}_{H_E}&=\inf\{\norm{f}_{\HH}:f\in\HH_+,\ f|_E=u\},\\
\norm{u}_{K_E}&=\inf\{\norm{f}_{\HH}:f\in\HH_-,\ f|_E=u\}.
\end{split}
\end{equation}
Both spaces contain every function on $E$. Indeed, prescribe the values $u(p)$ and $\pm u(p)$ at $p$ and $-p$, respectively. Smooth interpolation on $E\cup(-E)$, followed by the corresponding symmetrisation, gives an even or odd function with restriction $u$; see the bump-function construction in \cite[Chapter~2]{Lee}. Thus, $H_E$ and $K_E$ are the same vector space of functions, equipped with different Hilbert norms.

The restriction maps have closed kernels. By the Hilbert space projection theorem \cite[Theorems~I.2.5--I.2.7]{Conway}, each quotient has a linear isometric lift obtained by choosing the representative of least norm. In particular, there is an isometry
\[
L_E:K_E\longrightarrow\HH_-,\qquad (L_Eu)|_E=u.
\]
The constant function on $E$ has $H_E$-norm at most one.

For $p\in E$, evaluation on either quotient has the same norm as on the corresponding full space, by the isometric identification of the dual of a quotient with the annihilator of its kernel \cite[Theorem~III.10.2]{Conway}. On the $q$-fold sums, write
\[
e_p^H:H_E^{(q)}\longrightarrow\CC^q,
\qquad e_p^K:K_E^{(q)}\longrightarrow\CC^q
\]
for the componentwise evaluations. By \eqref{eq:evalnorms} and the Riesz representation theorem,
\begin{equation}\label{eq:coisometry}
e_p^H(e_p^H)^*=a^2I_q,
\qquad e_p^K(e_p^K)^*=b^2I_q.
\end{equation}

For any scalar function $\lambda:E\to\CC$, let $A_E(\lambda)$ and $B_E(\lambda)$ be multiplication by $\lambda$ on $H_E$ and $K_E$:
\begin{equation}\label{eq:multiplication}
A_E(\lambda)u=\lambda u,\qquad B_E(\lambda)u=\lambda u.
\end{equation}
If $\lambda$ is injective, both operators are diagonalizable with simple eigenvalues $\lambda(p)$. The eigenspace at $\lambda(p)$ consists of the functions supported at $p$. Pointwise bounds on $\lambda$ alone need not bound these operator norms uniformly as $E$ varies.

\subsection{Odd multiplicities}

We show that the norms of a similarity and its inverse cannot both stay bounded as the finite sets and their antipodes fill the sphere.

We call $E\cup(-E)$ an $\eta$-net if each point of the sphere is at distance less than $\eta$ from this set. We shall use the fact that every continuous odd function $\SSS\to\CC$ has a zero. This is the Borsuk--Ulam theorem applied to the real and imaginary parts; see \cite[Theorem~1.10, p.~32]{Hatcher}.

\begin{proposition}\label{prop:odd}
For every $M\ge1$, there is $\eta>0$ with the following property. Let $E$ contain no antipodal pair, suppose that $E\cup(-E)$ is an $\eta$-net, and let $\lambda:E\to\CC$ be injective. For every odd positive integer $q$, each invertible operator
\[
Y:H_E^{(q)}\longrightarrow K_E^{(q)},\qquad
Y A_E(\lambda)^{(q)}=B_E(\lambda)^{(q)}Y,
\]
satisfies
\begin{equation}\label{eq:largecondition}
\norm{Y}\norm{Y^{-1}}>M.
\end{equation}
The choice of $\eta$ is independent of $q$, of $E$, and of the values of $\lambda$.
\end{proposition}

The constant function and the evaluation operators give the required estimate.

\begin{proof}
Put $\gamma=b/a$. By \eqref{eq:evalcontinuous}, choose $\eta>0$ such that
\begin{equation}\label{eq:smallvariation}
|x-y|<\eta\quad\Longrightarrow\quad
\norm{e_x^- -e_y^-}<\frac{\gamma}{M}.
\end{equation}
Suppose that such a $Y$ has $\norm{Y}\norm{Y^{-1}}\le M$. After a positive scalar rescaling, we may assume
\begin{equation}\label{eq:balanced}
\norm{Y}=\norm{Y^{-1}}\le\sqrt M.
\end{equation}
Since the eigenvalues in \eqref{eq:multiplication} are distinct, the intertwining equation makes $Y$ act pointwise:
\begin{equation}\label{eq:pointwise}
(Yu)(p)=G_pu(p),\qquad G_p\in GL_q(\CC),\quad p\in E.
\end{equation}
In particular, $e_p^HY^{-1}=G_p^{-1}e_p^K$. Equation~\eqref{eq:coisometry} gives
\[
b\norm{G_p^{-1}}=\norm{G_p^{-1}e_p^K}
\le a\norm{Y^{-1}}\le a\sqrt M.
\]
Consequently,
\begin{equation}\label{eq:pointwiselower}
\norm{G_p\xi}\ge\frac{\gamma}{\sqrt M}\norm{\xi}
\qquad(\xi\in\CC^q).
\end{equation}

Let $C_E:\CC^q\to H_E^{(q)}$ send $\xi$ to the constant function with value $\xi$. Since $\norm{C_E}\le1$, the operator $U=L_E^{(q)}YC_E$ maps $\CC^q$ into $\HH_-^{(q)}$ and has norm at most $\sqrt M$. Define the matrix-valued function $\Phi$ by
\begin{equation}\label{eq:Phi}
\Phi(x)\xi=(U\xi)(x).
\end{equation}
It is continuous and odd, and $\Phi(p)=G_p$ for $p\in E$. Under the Hilbert space identification $\HH_-^{(q)}=\HH_-\otimes\CC^q$, componentwise evaluation is $e_x^-\otimes I_q$. Hence
\[
\norm{(e_x^--e_y^-)\otimes I_q}=\norm{e_x^--e_y^-}.
\]
It follows that
\begin{equation}\label{eq:variation}
\norm{\Phi(x)-\Phi(y)}\le\sqrt M\,\norm{e_x^- -e_y^-}.
\end{equation}
Since $q$ is odd, $\det\Phi$ is a continuous odd complex-valued function. Borsuk--Ulam gives a point $z$ where $\Phi(z)$ is singular. Choose a unit vector $\xi\in\ker\Phi(z)$ and a point $p\in E\cup(-E)$ with $|p-z|<\eta$. At this point, $\Phi(p)$ is $G_x$ or $-G_x$ for some $x\in E$. Equations \eqref{eq:pointwiselower}, \eqref{eq:smallvariation}, and~\eqref{eq:variation} give the contradiction
\[
\frac{\gamma}{\sqrt M}
\le\norm{\Phi(p)\xi}
=\norm{(\Phi(p)-\Phi(z))\xi}
<\sqrt M\,\frac{\gamma}{M}.\qedhere
\]
\end{proof}

\subsection{An invertible block intertwiner}

The intertwining equation for primitive square roots requires equal diagonal blocks; see \cite[Lemma~4.9]{MRTZ}. We choose those blocks, and the off-diagonal entries, so that the inverse can be written explicitly.

\begin{lemma}\label{lem:blockidentity}
Let $Z,W,V$ be commuting bounded operators on a Hilbert space such that
\[
Z^2+WV=I,
\]
and suppose that $R=I+WV$ is invertible. Set
\[
D=-V^2R^{-1},\qquad
F=\begin{pmatrix}Z&W\\DW&Z\end{pmatrix}.
\]
Then,
\begin{equation}\label{eq:blockinverse}
F^{-1}=\begin{pmatrix}RZ&-RW\\VWV&RZ\end{pmatrix},
\qquad FJ_2(D)=J_2(D)F.
\end{equation}
If a bounded involution $\Gamma$ anticommutes with $Z,W,V$, then $D$ commutes with $\Gamma$, while $F$ and $F^{-1}$ interchange the two eigenspaces of $\Gamma\oplus\Gamma$.
\end{lemma}

\begin{proof}
All the operators in the block entries commute. The relation $Z^2=I-WV$ gives
\[
R(Z^2-DW^2)=(I+WV)(I-WV)+(WV)^2=I.
\]
Multiplying the two block matrices in either order gives the identity, with lower-left inverse entry $-RDW=V^2W=VWV$. Comparing the products with $J_2(D)$ proves the intertwining identity. If $\Gamma$ anticommutes with $Z,W,V$, then $R,R^{-1},D$ commute with $\Gamma$, and each entry of $F$ and $F^{-1}$ anticommutes with $\Gamma$. This proves the last assertion.
\end{proof}

We realise these operators on $\HH$. For a smooth scalar function $g$, let $M_g$ denote multiplication by $g$, and put
\begin{equation}\label{eq:coordinateoperators}
Z=M_{x_3},\qquad W=M_{x_1+ix_2},\qquad V=M_{x_1-ix_2}.
\end{equation}
They are bounded, commute, and satisfy $Z^2+WV=I$ by the equation of the sphere. All three anticommute with $\Gamma$. Here $V$ is multiplication by the conjugate coordinate, not the Hilbert space adjoint of $W$.

The operator $R$ is multiplication by $1+x_1^2+x_2^2$. Its reciprocal is smooth, so $R^{-1}$ is a bounded multiplier. This argument does not require $R$ to be self-adjoint for the Sobolev inner product. The operator $D$ is multiplication by the smooth even function
\begin{equation}\label{eq:dsymbol}
d(x)=-\frac{(x_1-ix_2)^2}{1+x_1^2+x_2^2}.
\end{equation}
Lemma~\ref{lem:blockidentity} therefore gives a bounded invertible operator
\begin{equation}\label{eq:fullintertwiner}
F:\HH_+^{(2)}\longrightarrow\HH_-^{(2)}.
\end{equation}
Both $F$ and its inverse have smooth multiplication operators as entries. They are fixed independently of the finite set $E$.

We now pass this operator and its inverse to the finite-dimensional quotients.

\begin{proposition}\label{prop:quotientintertwiner}
There is $C\ge1$ such that, for each finite nonempty $E$ containing no antipodal pair, an invertible operator
\[
X_E:H_E^{(2)}\longrightarrow K_E^{(2)}
\]
satisfies $\norm{X_E},\norm{X_E^{-1}}\le C$ and
\begin{equation}\label{eq:allmultipliers}
X_EA_E(\lambda)^{(2)}=B_E(\lambda)^{(2)}X_E
\end{equation}
for every scalar function $\lambda$ on $E$. Moreover,
\begin{equation}\label{eq:quotientroots}
X_EJ_2(A_E(d|_E))=J_2(B_E(d|_E))X_E,
\end{equation}
and
\begin{equation}\label{eq:quotientDnorms}
\norm{A_E(d|_E)},\norm{B_E(d|_E)}\le\norm{D}.
\end{equation}
\end{proposition}

\begin{proof}
Multiplication preserves vanishing on $E$. Thus, $F$ maps the restriction kernel in $\HH_+^{(2)}$ into that in $\HH_-^{(2)}$, and $F^{-1}$ maps the latter kernel back into the former. The two inclusions show that $F$ maps the first kernel onto the second. Hence, the induced quotient maps are inverses. By the quotient norms,
\[
\norm{X_E}\le\norm{F},\qquad
\norm{X_E^{-1}}\le\norm{F^{-1}}.
\]
For example, $\norm{X_Eu}$ is at most $\norm{Ff}$ for every lift $f$ of $u$; taking the infimum over lifts gives the first bound. The inverse is treated in the same way.

At each $p\in E$, the matrix defining $X_E$ commutes with $\lambda(p)I_2$. This proves \eqref{eq:allmultipliers}, even when $\lambda$ is specified only on $E$. The identity with $J_2(D)$ in \eqref{eq:blockinverse} passes to the quotients and gives \eqref{eq:quotientroots}. Finally, $D$ preserves both parity spaces and both restriction kernels. Its quotient operators have norm at most $\norm{D}$, proving \eqref{eq:quotientDnorms}.
\end{proof}

The distinction between the two intertwining equations will be used below. Equation~\eqref{eq:allmultipliers} allows any finite choice of eigenvalues. Equation~\eqref{eq:quotientroots} requires the particular multiplier $D$.

\section{Counterexamples on Hilbert space}\label{sec:examples}

We first choose the finite sets and recall how disjoint spectra control the intertwiners of an orthogonal sum. Both constructions will then use the operators from Section~\ref{sec:quotients}.

\subsection{Finite sets and spectral separation}

The finite sets must become dense after adjoining their antipodes. We also choose them so that no eigenvalue occurs in two different summands.

Let
\[
U=\{x\in\SSS:0<x_3<1\}.
\]
Choose finite nonempty sets $E_m\subset U$ such that
\begin{equation}\label{eq:nets}
\text{the covering radii of }E_m\cup(-E_m)\text{ tend to zero,}
\end{equation}
and all third coordinates of all selected points are distinct. Such sets can be chosen inductively. For a prescribed radius, compactness of the sphere and density of $U\cup(-U)$ give a finite antipodal net with smaller radius. Perturb its points within $U$ to avoid the third coordinates chosen earlier and to make the new coordinates distinct. Only finitely many values are excluded, and every nonempty open subset of $U$ contains an interval of possible third coordinates. The perturbations can be made small enough to retain the prescribed covering radius.

Suppose that $A_m$ and $B_m$ act on finite-dimensional Hilbert spaces, their norms are uniformly bounded, and
\begin{equation}\label{eq:separation}
\sigma(A_m)=\sigma(B_m)=\Sigma_m,
\qquad \Sigma_m\cap\Sigma_n=\varnothing\quad(m\ne n).
\end{equation}
If $Y(\bigoplus_m A_m)=(\bigoplus_m B_m)Y$, then its blocks satisfy
\[
Y_{nm}A_m=B_nY_{nm}.
\]
For $n\ne m$, the Sylvester--Rosenblum theorem \cite{Rosenblum} gives $Y_{nm}=0$. Thus, $Y$ is block diagonal. If $Y$ is invertible, its diagonal blocks are invertible and
\begin{equation}\label{eq:blocknorms}
\norm{Y_m}\le\norm{Y},\qquad
\norm{Y_m^{-1}}\le\norm{Y^{-1}}.
\end{equation}
The same conclusion holds after a common finite amplification. This is the orthogonal-sum version of the spectral separation used in \cite[Lemma~2.16]{MRTZ}. Only disjointness of the finite spectra is needed; their union need not have isolated points.

\subsection{Similar primitive square roots}

We are now ready to prove Theorem~\ref{thm:bounded}. The multiplier $D$ gives uniformly bounded finite-dimensional operators, and the maps $X_E$ give a similarity of their primitive square roots.

\begin{proof}[Proof of Theorem~\ref{thm:bounded}]
On $H_{E_m}$ and $K_{E_m}$, respectively, set
\[
A_m=A_{E_m}(d|_{E_m}),\qquad B_m=B_{E_m}(d|_{E_m}),
\]
and define
\begin{equation}\label{eq:boundedpair}
T_0=\bigoplus_m A_m\quad\text{on }\bigoplus_m H_{E_m},
\qquad
S_0=\bigoplus_m B_m\quad\text{on }\bigoplus_m K_{E_m}.
\end{equation}
Both sums are bounded by \eqref{eq:quotientDnorms}. Proposition~\ref{prop:quotientintertwiner} also shows that
\[
X=\bigoplus_m X_{E_m},\qquad
X^{-1}=\bigoplus_m X_{E_m}^{-1}
\]
are bounded. Reordering the orthogonal sums and applying \eqref{eq:quotientroots} gives
\begin{equation}\label{eq:globalroot}
XJ_2(T_0)=J_2(S_0)X.
\end{equation}
Taking squares gives similarity of the doubles. Repetition gives similarity at every even multiplicity.

For $x\in U$, the function in \eqref{eq:dsymbol} satisfies
\begin{equation}\label{eq:modulus}
|d(x)|=\frac{1-x_3^2}{2-x_3^2}.
\end{equation}
This is strictly decreasing as $x_3$ increases from zero to one. Hence, the values $d(x)$ at all selected points are nonzero and pairwise distinct, including values from different sets. The finite blocks in \eqref{eq:boundedpair} therefore satisfy \eqref{eq:separation}.

Suppose that an invertible $Y$ intertwines an odd amplification of $T_0$ and $S_0$. Spectral separation gives invertible restrictions $Y_m$ satisfying \eqref{eq:blocknorms}. Choose $M\ge1$ with $M\ge\norm{Y}\norm{Y^{-1}}$. By \eqref{eq:nets}, every sufficiently large $E_m$ satisfies the net condition in Proposition~\ref{prop:odd}. Hence,
\[
\norm{Y_m}\norm{Y_m^{-1}}>M,
\]
contrary to \eqref{eq:blocknorms}. We have proved
\begin{equation}\label{eq:boundedparity}
T_0^{(q)}\sim S_0^{(q)}\quad\Longleftrightarrow\quad q\text{ is even}.
\end{equation}

The point spectrum of either operator is exactly $\bigcup_m d(E_m)$. Indeed, the components of an eigenvector must solve the corresponding finite-block eigenvalue equations. Since all the eigenvalues are distinct, every eigenspace is one-dimensional. The eigenvectors span each finite block, so their linear span is dense in each direct sum.

Neither operator is compact. Fix a point on the equator and choose points in $E_m\cup(-E_m)$ converging to it. Their third coordinates tend to zero, so the corresponding eigenvalue moduli tend to $1/2$ by \eqref{eq:modulus}. For either operator, choose unit eigenvectors in these distinct finite blocks. They form an orthonormal sequence whose images do not tend to zero. The compact-operator criterion in \cite[Chapter~II, \S4]{Conway} excludes compactness.

To obtain the assertion about all primitive roots, we use two identities from \cite{MRTZ}: similarity of $A$ and $B$ implies similarity of $J_k(A)$ and $J_k(B)$ by Lemma~2.4(c), and
\[
J_k(J_n(A))\simeq_u J_{kn}(A)
\]
by Proposition~2.5(a). Applying these facts to \eqref{eq:globalroot} proves similarity of the primitive roots at every even order. At an odd order $n$, such a similarity would give $T_0^{(n)}\sim S_0^{(n)}$ after taking $n$th powers, contrary to \eqref{eq:boundedparity}.

For $t>0$, let $L_t=\diag(t^{-1/2}I,I)$. The scalar-rescaling identity recorded in \cite[Remark~5.3]{MRTZ} is
\begin{equation}\label{eq:rescale}
J_2(tA)=t^{1/2}L_tJ_2(A)L_t^{-1}.
\end{equation}
Choose $t>0$ with $t\max\{\norm{T_0},\norm{S_0}\}<\varepsilon$, and put $T=tT_0$, $S=tS_0$. The intertwining equations for finite amplifications are unchanged by a common nonzero scalar. Equation~\eqref{eq:rescale} preserves the similarity of the primitive square roots. The preceding even-order argument and the odd-power argument then give \eqref{eq:rootparity} for this scaled pair. The assertions about point spectra, eigenspaces, and noncompactness follow from those proved above. Finally, identify the two separable infinite-dimensional Hilbert spaces unitarily with $H$.
\end{proof}

\subsection{Compact operators with prescribed singular value bounds}

For the compact pair, we keep the quotient spaces and the similarities of the doubles, but choose positive eigenvalues separately on each finite set. Small scalar multiples will give the singular-value bounds in Theorem~\ref{thm:compact}.

\begin{proof}[Proof of Theorem~\ref{thm:compact}]
Choose any injective function
\[
\lambda_m:E_m\longrightarrow(0,\infty)
\]
on each finite set, and now put
\[
A_m=A_{E_m}(\lambda_m),\qquad B_m=B_{E_m}(\lambda_m).
\]
Their norms are finite, although they may grow with $m$. Proposition~\ref{prop:quotientintertwiner} gives the same uniformly bounded similarities of their doubles, while Proposition~\ref{prop:odd} still applies at every odd multiplicity.

Write $n_m=|E_m|$, $N_0=0$, and $N_m=n_1+\cdots+n_m$. Choose positive numbers $t_m$ inductively so that
\begin{equation}\label{eq:weights}
t_m\max\{\norm{A_m},\norm{B_m}\}<a_{N_m}
\end{equation}
and the finite sets $t_m\lambda_m(E_m)$ are pairwise disjoint. At each step, \eqref{eq:weights} leaves a nonempty interval of choices, whereas spectral disjointness excludes only finitely many ratios of previous eigenvalues to the values of $\lambda_m$. Hence, both requirements can be satisfied.

Define
\begin{equation}\label{eq:compactpair}
K=\bigoplus_m t_mA_m,\qquad L=\bigoplus_m t_mB_m.
\end{equation}
By \eqref{eq:weights}, the tails tend to zero in operator norm, so $K$ and $L$ are compact. For $N_{m-1}<j\le N_m$, the first $m-1$ blocks have total rank less than $j$, and the remaining tail has norm at most
\[
\sup_{k\ge m}a_{N_k}=a_{N_m}\le a_j.
\]
The singular-value approximation principle \cite[Chapter~1]{Simon},
\[
s_j(R)=\inf\{\norm{R-F}:\rank F<j\},
\]
therefore gives \eqref{eq:singularbounds} for both operators.

The maps $X_{E_m}$ and their inverses are uniformly bounded. Their direct sums are bounded inverses of one another, and \eqref{eq:allmultipliers}, applied to $t_m\lambda_m$, gives $K^{(2)}\sim L^{(2)}$. Repetition gives all even multiplicities. If an invertible $Y$ intertwined an odd amplification, the disjoint finite spectra would make it block diagonal. Its restrictions would satisfy \eqref{eq:blocknorms}, whereas Proposition~\ref{prop:odd} contradicts those bounds for sufficiently large $m$. This proves \eqref{eq:compactparity}.

Every finite block is invertible. Thus, $K,L$ are injective, and their ranges contain all finitely supported vectors, so both ranges are dense. Their eigenvectors span each finite block and have dense linear span. By the Riesz--Schauder theorem \cite[Chapter~VII, \S7]{Conway}, every nonzero spectral value of a compact operator is an eigenvalue. The block equations identify these values as
\[
\bigcup_m t_m\lambda_m(E_m).
\]
They are positive and distinct. Each belongs to a single diagonalizable finite block. The same block equations for powers of $K-\mu I$ and $L-\mu I$ show that the generalized eigenspace at each nonzero eigenvalue $\mu$ is one-dimensional. Hence, each such eigenvalue is algebraically simple. Unitary identifications put both operators on $H$.
\end{proof}

For $p>0$, a compact operator $R$ belongs to $\mathcal S_p$ when $\sum_j s_j(R)^p<\infty$; see \cite[Chapter~1]{Simon}. Choosing a sequence whose positive powers are all summable gives the following consequence.

\begin{corollary}\label{cor:schatten}
There are injective compact operators $K,L$ with dense range, belonging to every Schatten class, such that $K^{(2)}\sim L^{(2)}$ but $K\not\sim L$.
\end{corollary}

\begin{proof}
This follows from Theorem~\ref{thm:compact} by taking $a_j=2^{-j}$.
\end{proof}

For the compact pair, $J_2(K)\not\sim J_2(L)$ by \cite[Theorem~4.10]{MRTZ}. The two constructions use the same quotient spaces and the same doubled intertwiners. Changing the eigenvalues and choosing small weights preserves \eqref{eq:allmultipliers}, but not the additional identity \eqref{eq:quotientroots} required for primitive square roots.

\section{Matrix size and bounded products}\label{sec:products}

We first compare the similarities of the finite-dimensional pairs before and after amplification. We then use those pairs to decide cancellation in bounded products of matrices and in finite type~I von Neumann algebras.

The following corollary records the finite-dimensional pairs obtained from Propositions~\ref{prop:odd} and~\ref{prop:quotientintertwiner}.

\begin{corollary}\label{prop:finite-bad}
There is $C\ge1$ with the following property. For every $M\ge1$, there are invertible diagonalizable matrices $A,B$ of the same size, with distinct eigenvalues and $\norm{A},\norm{B}\le1$, and an invertible $X$ such that
\begin{equation}\label{eq:finite-root-bounds}
XJ_2(A)=J_2(B)X,\qquad \norm{X},\norm{X^{-1}}\le C,
\end{equation}
while every invertible $Y$ satisfying $YA^{(q)}=B^{(q)}Y$ for an odd positive integer $q$ has
\begin{equation}\label{eq:finite-large}
\norm{Y}\norm{Y^{-1}}>M.
\end{equation}
The matrix size may depend on $M$.
\end{corollary}

\begin{proof}
Take the spaces $H_{E_m},K_{E_m}$ and multiplication by $d|_{E_m}$. For sufficiently large $m$, Proposition~\ref{prop:odd} gives \eqref{eq:finite-large} simultaneously for all odd $q$. The eigenvalues are distinct and nonzero by \eqref{eq:modulus} and the choice of the sets. Proposition~\ref{prop:quotientintertwiner} gives operator norms at most $\norm{D}$ and an intertwiner of the primitive square roots whose norm and inverse norm are uniformly bounded.

Choose one fixed scalar $c>0$ such that $c\norm{D}\le1$, and multiply both operators by $c$. This leaves the intertwining equations for finite amplifications unchanged. By \eqref{eq:rescale}, the primitive-root intertwiner is conjugated by $L_c$. Its norm and inverse norm are multiplied by at most $\norm{L_c}\norm{L_c^{-1}}$, a fixed number. Thus, \eqref{eq:finite-root-bounds} holds with a constant independent of $m$ and $M$. Represent the two finite-dimensional Hilbert spaces in orthonormal bases to obtain the stated matrices.
\end{proof}

Taking squares in \eqref{eq:finite-root-bounds} gives $s(A^{(2)},B^{(2)})\le C^2$. Equation~\eqref{eq:finite-large} with $q=1$, together with Theorem~\ref{thm:uniform}, proves \eqref{eq:size-dependence} with $t_0=C^2$. Thus, the failure of a bound independent of matrix size occurs even for invertible diagonalizable contractions with simple spectra, and even when their primitive square roots have uniformly bounded similarities. Taking $M>C^2$ in Corollary~\ref{prop:finite-bad} also gives $s(A^{(2)},B^{(2)})<s(A,B)$, so the equality suggested in \cite[Remark~4.2]{JMSZ} does not hold in general.

Let $(d_j)_{j\ge1}$ be a sequence of positive integers. In the theorem below,
\begin{equation}\label{eq:product-algebra}
\MM=\prod_{j\ge1}M_{d_j}(\CC)
\end{equation}
denotes the von Neumann algebra of all sequences $(A_j)$ with $\sup_j\norm{A_j}<\infty$, equipped with the supremum norm. This is the bounded product, not the norm-closed direct sum. An element $(Y_j)$ is invertible in $\MM$ precisely when every $Y_j$ is invertible and both $\sup_j\norm{Y_j}$ and $\sup_j\norm{Y_j^{-1}}$ are finite.

The reduction for the full sequence of matrix sizes at multiplicity two is given in \cite[Corollary~4.6]{JMSZ}. We now decide cancellation for each sequence of sizes and include the primitive-square-root question.

\begin{theorem}\label{thm:products}
For the algebra in \eqref{eq:product-algebra}, the following are equivalent.
\begin{enumerate}[label=\textup{(\roman*)},leftmargin=*,itemsep=3pt]
\item $\sup_jd_j<\infty$.
\item For all $A,B\in\MM$, similarity of $A^{(2)}$ and $B^{(2)}$ in $M_2(\MM)$ implies similarity of $A$ and $B$ in $\MM$.
\item For all $A,B\in\MM$, similarity of $J_2(A)$ and $J_2(B)$ in $M_2(\MM)$ implies similarity of $A$ and $B$ in $\MM$.
\end{enumerate}
When these conditions hold, cancellation is valid for every finite amplification and every primitive root.
\end{theorem}

The fixed-size theorem proves one direction. For the converse, the matrices in Corollary~\ref{prop:finite-bad} fit into a subsequence of factors of increasing size.

\begin{proof}
Assume~(i), and choose $N$ with $d_j\le N$ for every $j$. Suppose that $XA^{(q)}=B^{(q)}X$ for an invertible $X\in M_q(\MM)$, and put $t=\norm{X}\norm{X^{-1}}$. Evaluation in the $j$th factor gives
\[
s(A_j^{(q)},B_j^{(q)})\le t.
\]
By Theorem~\ref{thm:uniform},
\[
s(A_j,B_j)\le K:=\max_{1\le n\le N}\theta_{n,q}(t)<\infty.
\]
For each $j$, choose $Y_jA_j=B_jY_j$ with condition number less than $K+1$, and balance $Y_j$ by a positive scalar. Then,
\[
\norm{Y_j},\norm{Y_j^{-1}}<\sqrt{K+1}\qquad(j\ge1).
\]
Thus, $Y=(Y_j)$ is invertible in $\MM$ and satisfies $YA=BY$. This proves~(ii) and cancellation at every finite multiplicity. Taking powers gives cancellation for every primitive root. In particular, (ii) implies~(iii) directly by squaring.

Suppose now that $(d_j)$ is unbounded. For each $m\ge1$, choose matrices $A_m,B_m\in M_{k_m}(\CC)$ and an intertwiner $X_m$ as in Corollary~\ref{prop:finite-bad}, with $M=m$. Choose distinct indices $j_m$ inductively so that $d_{j_m}\ge k_m$. In the $j_m$th factor set
\[
\widehat A_{j_m}=A_m\oplus0_{d_{j_m}-k_m},\qquad
\widehat B_{j_m}=B_m\oplus0_{d_{j_m}-k_m},
\]
and put $\widehat A_j=\widehat B_j=0$ in the other factors. These sequences define contractions $\widehat A,\widehat B\in\MM$.

A permutation of the direct summands identifies $J_2(\widehat A_{j_m})$ with
\[
J_2(A_m)\oplus J_2(0_{d_{j_m}-k_m}),
\]
and gives the same identification for $\widehat B_{j_m}$. The map $X_m\oplus I$ intertwines these two operators. Its norm and inverse norm are at most $C$, after increasing the fixed constant to at least one. In the other factors use the identity. The resulting sequence is an invertible element of $M_2(\MM)$, so
\begin{equation}\label{eq:product-root-counter}
J_2(\widehat A)\sim J_2(\widehat B)
\quad\text{in }M_2(\MM).
\end{equation}

If $\widehat A$ and $\widehat B$ were similar in $\MM$, their coordinate similarities would have uniformly bounded condition numbers. In the $j_m$th factor, invertibility of $A_m$ and $B_m$ separates their spectra from the additional zero block. The Sylvester--Rosenblum theorem therefore forces any such similarity to restrict to a similarity between $A_m$ and $B_m$, with no increase in its norm or inverse norm. By \eqref{eq:finite-large}, the condition number of this restriction is greater than $m$. This contradicts the uniform bound. Thus, \eqref{eq:product-root-counter} holds without similarity of $\widehat A$ and $\widehat B$, which disproves~(iii). We have shown that (iii) implies~(i), completing the proof.
\end{proof}

\subsection{Finite type I von Neumann algebras}
The centres of the homogeneous summands do not change the conclusion. We use the central decomposition of a finite type~I von Neumann algebra \cite[\S6.6]{KR}, in the form
\begin{equation}\label{eq:finiteI-decomposition}
\MM=\prod_{n\in J}M_n(\mathcal A_n),
\end{equation}
where $J\subseteq\mathbb N$ and each $\mathcal A_n$ is a nonzero abelian von Neumann algebra. The product in \eqref{eq:finiteI-decomposition} is again the bounded product. We exclude the zero algebra, so $J$ is nonempty.

\begin{theorem}\label{thm:finiteI}
Let $\MM$ be a nonzero finite type~I von Neumann algebra, with decomposition \eqref{eq:finiteI-decomposition}. The following conditions are equivalent.
\begin{enumerate}[label=\textup{(\roman*)},leftmargin=*,itemsep=3pt]
\item The set $J$ is bounded.
\item For all $A,B\in\MM$, similarity of $A^{(2)}$ and $B^{(2)}$ in $M_2(\MM)$ implies similarity of $A$ and $B$ in $\MM$.
\item For all $A,B\in\MM$, similarity of $J_2(A)$ and $J_2(B)$ in $M_2(\MM)$ implies similarity of $A$ and $B$ in $\MM$.
\end{enumerate}
When these conditions hold, cancellation is valid for every finite amplification and every primitive root.
\end{theorem}

To prove the converse, it is important not to infer nonsimilarity in $\MM$ merely from nonsimilarity in a subalgebra. The abelian centres provide homomorphisms back to the scalar matrix factors.

\begin{proof}
If $J$ is bounded, it is finite. Apply Theorem~\ref{thm:typeI} in each of the finitely many homogeneous summands. The resulting similarities and their inverses form bounded families, so they give a similarity in $\MM$. This proves cancellation for every finite amplification. Taking powers gives the assertion for primitive roots. As before, (ii) implies (iii) by squaring.

Suppose that $J$ is unbounded, and put
\[
\mathcal N=\prod_{n\in J}M_n(\CC).
\]
By Theorem~\ref{thm:products}, there are $A,B\in\mathcal N$ whose primitive square roots are similar in $M_2(\mathcal N)$ but which are not similar in $\mathcal N$. Let
\[
\iota:\mathcal N\longrightarrow\MM,
\qquad
\iota((A_n))=(A_n\otimes1_{\mathcal A_n})
\]
be the isometric unital inclusion. It carries the given similarity to one between $J_2(\iota(A))$ and $J_2(\iota(B))$.

For each $n$, choose a character $\chi_n$ on $\mathcal A_n$. Such a character exists because $\mathcal A_n$ is a nonzero unital abelian $C^*$-algebra; see \cite[\S I.3]{Davidson}. The characters are not assumed to be normal. Applying $\chi_n$ entrywise gives a unital $*$-homomorphism
\[
\pi:\MM\longrightarrow\mathcal N,
\qquad
\pi((T_n))=
\bigl((\operatorname{id}_{M_n}\otimes\chi_n)(T_n)\bigr)_{n\in J}.
\]
Each coordinate map is contractive, so $\pi$ is well defined on the bounded product, and $\pi\iota=\operatorname{id}_{\mathcal N}$. If an invertible $Y\in\MM$ intertwined $\iota(A)$ and $\iota(B)$, then $\pi(Y)$ would intertwine $A$ and $B$ and would have inverse $\pi(Y^{-1})$. This contradicts their nonsimilarity in $\mathcal N$. Hence, (iii) fails when $J$ is unbounded, proving the theorem.
\end{proof}

\end{document}